\documentclass{article}
\usepackage{amssymb}
\usepackage{amsmath}
\usepackage{amsthm}
\usepackage{amsbsy}
\usepackage{psfrag}
\usepackage{graphics}
\usepackage{graphicx}

\newtheorem{theorem}{Theorem}[section]

\theoremstyle{definition}

\theoremstyle{remark}

\numberwithin{equation}{section}
\newcommand{\bbP}{\mathbb{P}}

\newcommand{\x}{\textbf{x}}
\renewcommand{\v}{\textbf{v}}
\renewcommand{\u}{\textbf{u}}
\newcommand{\n}{\textbf{n}}
\usepackage{graphicx}
\usepackage{psfrag}
\usepackage{caption2}
\begin{document}

\title{Optimal lower bounds on the local stress inside random thermoelastic composites}

\author{Yue Chen$^{\mbox{\tiny 1}}$
\and Robert Lipton $^{\mbox{\tiny 2}}$\thanks{This work is partially supported by  grants: NSF DMS-0807265 and AFOSR FA9550-05-0008.}}
\date{}

\maketitle

\baselineskip=0.9\normalbaselineskip
\vspace{-3pt}

\begin{center}
{\footnotesize 
$^{\mbox{\tiny\rm 1}}$ Department of Mathematics 
Louisiana State University
Baton Rouge, LA\\ 70803, USA. email: chenyue\symbol{'100}math.lsu.edu\\[3pt]
$^{\mbox{\tiny\rm 2}}$Department of Mathematics, Louisiana State University,
Baton Rouge, LA\\ 70803, USA. email: lipton\symbol{'100}math.lsu.edu\\[3pt]}
\end{center}

\maketitle

\begin{abstract}
A methodology is presented for bounding all higher moments of the local hydrostatic stress field inside random two phase linear thermoelastic media undergoing macroscopic thermomechanical loading. The method also provides a lower bound on the maximum local stress. Explicit formulas for the optimal lower bounds are found that are expressed in terms of the applied macroscopic thermal and mechanical loading, coefficients of thermal expansion, elastic properties, and volume fractions. These bounds provide a means to measure load transfer across length scales relating the excursions of the local fields to the applied loads and the thermal stresses inside each phase. These bounds are shown to be the best possible in that they are attained by the Hashin-Shtrikman coated sphere assemblage.
\end{abstract}

\section{Introduction}
\label{Introduction}

Over the last century major strides have been made in the characterization of effective constitutive laws relating average fluxes to average gradients inside random heterogeneous media  see for example \cite{Hash, Maxwell04, Willis, Milton, NNH, Rayleigh}. However much less is known about the point wise behavior of local fluxes and gradients fields inside random media. While it is true that efficient numerical methods capable of resolving local fields are available for prescribed microstructures, it is also true that for most applications only a partial statistical description of the microstructure is available.  Thus for these cases one must resort to bounds or approximations for the local fields that are based upon the available statistical descriptors of the microgeometry and the applied macroscopic loading. 
Bounds are useful as they provide a means to quantitatively assess load transfer across length scales relating the excursions of the local fields to applied macroscopic loads. Moreover, they provide explicit criteria on the applied loads that are necessary for failure initiation inside statistically defined heterogeneous media \cite{AlaliLipton}. In this paper we develop lower bounds on local field properties for statistically defined two phase microstructures when only the volume fraction of each  phase is known.  
Here the focus is on lower bounds since volume constraints alone do not preclude the existence of microstructures with rough interfaces for which the $L^p$ norms of local fields are divergent see \cite{SeminalMilton}, \cite{faraco}, and also \cite{leonetinesi}.

We present a methodology for bounding the  $L^p$  norms, $2\leq p\leq\infty$, of the local hydrostatic stress field inside random media made up of two thermoelastic materials. The method is used to obtain new optimal lower bounds that are given by explicit formulas expressed in terms of the applied thermal and mechanical loading, coefficients  of thermal expansion, elastic properties, and volume fractions. We show that these bounds are the best possible in that they are attained by the local fields inside the coated sphere assemblage originally introduced in \cite{hscoated}. 
It has been  known since 1963 that the coated spheres microstructure exhibits extreme effective elastic properties \cite{hashin}. However it was discovered only recently in \cite{Lipstress}, \cite{Lipstrain} that this geometry supports extreme local fields that minimize the maximum local hydrostatic field over all two phase elastic mixtures in fixed volume fractions. More recently several scenarios are identified for which, in the absence of thermal stresses, these microstructures attain lower bounds on the total local stress field inside each material when the composite is subjected to mechanical loading see, \cite{AlaliLipton}.

In this paper we consider mixtures of two thermoelastic materials with shear and bulk moduli specified by $\mu_1$, $k_1$, $\mu_2$, $k_2$ and coefficients of thermal expansion given by $h_1$ and $h_2$. New lower bounds are presented for elastically well ordered phases for which $k_1>k_2$ and $\mu_1>\mu_2$ as well as for non well ordered phases such that $k_2>k_1$ and $\mu_1>\mu_2$. For each of these cases we consider both macroscopic mechanical and thermal loads and present bounds that hold for $h_1>h_2$ and $h_2>h_1$. The set of bounds and optimal microstructures for the well ordered case  are listed in Section 3 and optimal lower bounds for the non well ordered case are listed in Section 4. 
The methodology for deriving the bounds is presented in Section 5.

The optimal bounds and the associated coated sphere microstructures given in Sections 3 and 4  show  that there are combinations of applied stress and imposed temperature change for which the local hydrostatic stress inside the connected phase of the coated sphere assemblage vanishes identically. Other loading combinations are seen to cause the stress inside the included phase of the coated sphere assemblage to vanish identically.
Thus for these cases the applied hydrostatic stress is converted into a pure local shear stress inside a preselected phase. 

Recent related work provides optimal lower bounds on  local fields in the absence of thermal loads.
The work presented in \cite{AlaliLipton} provides new optimal lower bounds on both the local shear stress and the local hydrostatic component of stress for random media subjected to a series of progressively more general applied macroscopic stresses. These bounds are explicit and given in terms of volume fractions, elastic constants of each phase, and the applied macroscopic stress.
Earlier work considers random two phase elastic composites subject to imposed macroscopic hydrostatic stress and strain  see, \cite{Lipstrain} and \cite{Lipstress}, as well as dielectric composites subjected applied constant electric fields see, \cite{Lipelect}. Those efforts deliver optimal lower bounds on the $L^p$ norms for the hydrostatic components of local stress and strain fields as well as the magnitude of the local electric field for all $p$ in the range $2\leq p \leq \infty$. Other work examines the stress field around a single simply connected stiff elastic inclusion subjected to a remote constant stress at infinity \cite{Wheeler} and provides optimal lower bounds for the supremum of the
maximum principal stress. The work presented in \cite{GrabovskyandKohn} provides an optimal lower bound on the supremum of the maximum principal stress for two-dimensional periodic composites consisting of a single simply connected elastically stiff inclusion inside the period cell. The recent work of  \cite{He} builds on the earlier work of \cite{Lipstrain, Lipstress} and develops new lower bounds on the
$L^p$ norm of the local stress and strain fields inside statistically isotropic two-phase elastic composites. However to date those bounds have been shown to be optimal for $p=2$ see, \cite{He}. Their optimality for $p>2$ remains to be seen. Optimal upper and lower bounds on the $L^2$ norm of local gradient fields
are established using integral representation formulas in \cite{lipsiama}.

We conclude by providing the notation and summation conventions used in this article. Contractions of
stress or strain fields $\sigma$ and $\epsilon$ are defined by
$\sigma:\epsilon=\sigma_{ij}\epsilon_{ij}$ and $|\sigma|^2=\sigma:\sigma$, where
repeated indices indicate summation. Products of fourth order
tensors $C$ and strain tensors $\epsilon$ are written as
$C\epsilon$ and are given by $[C\epsilon]_{ij}=C_{ijkl}\epsilon_{kl}$; and
products of stresses $\sigma$ with vectors ${\bf v}$ are
given by $[\sigma{\bf v}]_i=\sigma_{ij}v_j$. The fourth order identity
map on the space of stresses or strains is denoted by $\mathbf I$
and ${\mathbf
I}_{ijkl}=1/2(\delta_{ik}\delta_{jl}+\delta_{il}\delta_{jk})$. The
projection onto the hydrostatic part of $\sigma(\x)$ is denoted by
$\bbP^H$ and is given explicitly by
\begin{eqnarray}
\bbP^H_{ijkl}=\frac{1}{d}\delta_{ij}\delta_{kl}&\,\hbox{and}\,&\bbP^H\psi(\x)=\frac{tr\,\sigma(\x)}{d}I.
\label{proj}
\end{eqnarray}
The projection onto the deviatoric part of $\sigma(\x)$ is denoted by
$\bbP^D$ and ${\mathbf I}=\bbP^H+\bbP^D$ with $\bbP^D\bbP^H=\bbP^H\bbP^D=0$. The tensor product between two vectors $\u$ and $\v$ is the matrix $\u\otimes\v$ with elements$[\u\otimes\v]_{ij}=u_i v_j$. 
Last we denote the basis  for the space of constant $3\times 3$ symmetric strain tensors by $\bar{\epsilon}^{kl}$ where $\bar{\epsilon}_{mn}^{kl}=\delta_{mk}\delta_{nl}$.

\bigskip
\section{Stress and strain fields inside heterogeneous thermoelastic\\ media with imposed macroscopic loading}
\label{Prestress Problem}

Several distinct physical processes can generate prestress within heterogeneous media.  In many cases it is generated  by a mismatch between the coefficients of thermal expansion of the component materials. To fix ideas we present the
physical model associated with this situation. The tensors of thermal expansion inside each phase are given by $\lambda_1=h_1 I$ and $\lambda_2=h_2 I$ where $I$ is the $3\times 3$ identity.  The elastic properties of each component are specified by the elasticity tensors $C_1$ and $C_2$ respectively. In this treatment we consider heterogeneous elastically isotropic materials and
the elasticity tensors of materials one and two are specified by 
\begin{eqnarray}
C_i=3k_i\mathbb{P}^H+2\mu_i\mathbb{P}^D, \hbox{ $i=1,2.$}
\label{elastconst}
\end{eqnarray} 
Without loss of generality we adopt the convention
\begin{eqnarray}
\mu_1>\mu_2.
\label{order}
\end{eqnarray}

The elastic displacement inside the composite is denoted by $\textbf{u}$ and the associated strain tensor is denoted by $\epsilon (\textbf{u})$. The position dependent elastic tensor and thermal expansion tensor for the heterogeneous medium are denoted by $C(\x)$ and $\lambda(\x)$ respectively where $\textbf{x}$ denotes a point inside the medium. The domain containing the composite is given by a cube $Q$ of unit side length. Here is supposed that $Q$ is the period cell for an infinite elastic medium. In what follows the integral of a quantity $q$ over the unit cube $Q$ is denoted by $\langle q \rangle$. 

A constant macroscopic stress $\overline{\sigma}$ and uniform change in temperature $\Delta T$ is imposed upon the heterogeneous material. The local stress inside the heterogeneous medium is expressed as the sum of a periodic mean zero fluctuation $\hat{\sigma}$
and $\overline{\sigma}$, i.e., $\sigma(\x)=\overline{\sigma}+\hat\sigma(\x)$ with $\langle\hat{\sigma}\rangle=0$. Elastic equilibrium inside each phase is given by:
\begin{eqnarray}
div \sigma =0.
\label{equlib}
\end{eqnarray}
The local elastic strain $\epsilon(\u)$ is related to the local stress through the constitutive law
\begin{eqnarray}
 \sigma(\textbf{x})=C(\textbf{x})(\epsilon (\textbf{u}(\textbf{x}))-\lambda(\textbf{x})\Delta T),
 \label{constitutive}
\end{eqnarray}
and the local elastic field is written in the form
\begin{eqnarray}
\epsilon(\u)=\overline{\epsilon}+\epsilon(\u^{per})
\label{epsilon}
\end{eqnarray}
where $\u^{per}$ is $Q$ periodic taking the same values on opposite sides of the period cell
and $\langle\epsilon(\u^{per})\rangle=0$.
Perfect contact between the component materials is  assumed, thus
both the displacement ${\u}$ and  traction $\sigma\n$ are
continuous across the two phase interface, i.e.,
\begin{eqnarray}
\label{equlibbdry1}
{\u}_{|_{\scriptscriptstyle{1}}}&=&{\u}_{|_{\scriptscriptstyle{2}}},\,\,\\
\label{equlibbdry2}
\sigma_{|_{\scriptscriptstyle{1}}} \n&=&\sigma_{|_{\scriptscriptstyle{2}}} \n.
\end{eqnarray}
Here the subscripts indicate the side of the interface
that the displacement and traction fields are evaluated on and  $\n$ denotes the normal vector to the interface pointing from material one into material two. 

The effective ``macroscopic''  constitutive law for the heterogeneous medium is given by the constant effective elasticity tensor $C^e$  and effective thermal stress tensor $H^e$ that provide the linear relation  between the imposed macroscopic stress $\bar{\sigma}$, uniform change in temperature $\Delta T$, and the average strain $\bar{\epsilon}$ given by \cite{Milton},
\begin{eqnarray}
\bar{\sigma}=C^e\bar{\epsilon}+H^e \Delta T.
\label{effelast}
\end{eqnarray}
Here the components of $C^e$ are given by
\begin{eqnarray}
C^e_{ijkl}=\langle C_{ijmn}(x)(\epsilon(\varphi^{kl})_{mn}+\bar{\epsilon}^{kl}_{mn})\rangle,
\label{effectelast}
\end{eqnarray}
where the fields $\varphi^{ij}$ are the periodic solutions of
 \begin{eqnarray}
div (C(x)(\epsilon(\varphi^{kl})+\bar{\epsilon}^{kl}))= 0  \quad\text{inside each phase }, 
 \label{basis}
 \end{eqnarray}
 with the appropriate traction and continuity conditions along the two phase interface given by
 \begin{eqnarray}
\label{equlibbdry1e}
{\varphi^{kl}}_{|_{\scriptscriptstyle{1}}}&=&{\varphi^{kl}}_{|_{\scriptscriptstyle{2}}},\,\,\\
\label{equlibbdry2e}
C_1\left(\epsilon(\varphi^{kl})+\bar{\epsilon}^{kl}\right)_{|_{\scriptscriptstyle{1}}} \n&=&C_2\left(\epsilon(\varphi^{kl})+\bar{\epsilon}^{kl}\right)_{|_{\scriptscriptstyle{2}}} \n.
\end{eqnarray}
The effective thermal stress tensor $H^e$ is given by 
\begin{eqnarray}
H^e=\langle C(x)(\epsilon(\varphi^p)-\lambda(\x)) \rangle.
\label{effectethstress}
\end{eqnarray}
Where $\varphi^p$ is the periodic solution of
\begin{eqnarray}
div (C(x)(\epsilon({\varphi}^p)-\lambda(\x)))= 0  \quad\text{inside each phase  },
\label{prestressequlib}
\end{eqnarray}
with the traction and continuity conditions along the two phase interface given by
\begin{eqnarray}
\label{equlibbdry1p}
{\varphi^p}_{|_{\scriptscriptstyle{1}}}&=&{\varphi^p}_{|_{\scriptscriptstyle{2}}},\,\,\\
\label{equlibbdry2p}
C_1\left(\epsilon(\varphi^p)-\lambda_1\right)_{|_{\scriptscriptstyle{1}}} \n&=&C_2\left(\epsilon(\varphi^p)-\lambda_2\right)_{|_{\scriptscriptstyle{2}}} \n.
\end{eqnarray}
From linearity it follows that the local fluctuating strain field is
given by the sum 
\begin{eqnarray}
\epsilon(\u^{per})=\epsilon(\varphi^{kl})\bar{\epsilon}^{kl}+\epsilon(\varphi^p)\Delta T.
\label{linear}
\end{eqnarray}
We write $\varphi^e=\varphi^{kl}\bar{\epsilon}^{kl}$ and from linearity one has the expression for $C^e\bar{\epsilon}$ given by
\begin{eqnarray}
\label{otherconstitut}
C^e\bar{\epsilon}=\langle C(\x)(\epsilon(\varphi^e)+\bar{\epsilon}\rangle.
\end{eqnarray}

In this article the imposed mechanical stress is given by a constant hydrostatic stress 
\begin{eqnarray}
\bar{\sigma}=\sigma_0 I,
\label{imposedstress}
\end{eqnarray}
where $\sigma_0$ can assume any value in $-\infty<\sigma_0<\infty$. 
We introduce a method for obtaining optimal lower bounds on the higher moments of the hydrostatic component
of the local stress inside the composite when it is subjected to an imposed hydrostatic load $\sigma_0 I$ and temperature change $\Delta T$. Here no restriction is placed on $\Delta T$.
The volume fractions of materials one and two are denoted by $\theta_1$ and $\theta_2$ and the average of a quantity $q$ over material one is denoted by $\langle q\rangle_1$ and over material two by $\langle q \rangle_2$.
In the following section we present optimal lower bounds on the following moments of the local hydrostatic stress $\bbP^H\sigma(\x)$ over the domain occupied by each material
given by
\begin{eqnarray}
\langle|\bbP^H\sigma(\x)|^p\rangle_1^{1/p} &\hbox{  and  }&\langle|\bbP^H\sigma(\x)|^p\rangle_2^{1/p},
\label{locstressmoment}
\end{eqnarray}
for $1< p\leq\infty$,
as well as for the maximum local hydrostatic stress over the whole composite domain
\begin{eqnarray}
\max_{\x \hbox{\tiny in Q}}\left\{|\bbP^H\sigma(\x)|\right\}.
\label{locstressmaxeverything}
\end{eqnarray}

It is pointed out that the case corresponding to $p=\infty$ in (\ref{locstressmoment}) corresponds to
lower bounds on the maximum local stress over each phase
\begin{eqnarray}
\max_{\x\hbox{ \tiny in material 1}}\left\{|\bbP^H\sigma(\x)|\right\}, && \max_{\x\hbox{ \tiny in material 2}}\left\{|\bbP^H\sigma(\x)|\right\},
\label{locstressmaxphase}
\end{eqnarray}

The lower bounds are given in terms of the volume fractions $\theta_1$ and $\theta_2$, as well as the bulk and shear moduli $k_1$, $k_2$, $\mu_1$, $\mu_2$, and the coefficients of thermal expansion $h_1$ and $h_2$. The lower bounds are described by the following characteristic combinations of these parameters given by:
\begin{eqnarray}
L_1&=& \frac{k_1(k_2+\frac{4}{3}\mu_2)}{k_1k_2+(k_1\theta_1+k_2\theta_2)\frac{4}{3}\mu_2},\label{L1}\\
L_2&=& \frac{k_2(k_1+\frac{4}{3}\mu_1)}{k_1k_2+(k_1\theta_1+k_2\theta_2)\frac{4}{3}\mu_1},\label{L2}\\
M_1&=& \frac{k_1(k_2+\frac{4}{3}\mu_1)}{k_1k_2+(k_1\theta_1+k_2\theta_2)\frac{4}{3}\mu_1},\label{M1}\\
M_2&=& \frac{k_2(k_1+\frac{4}{3}\mu_2)}{k_1k_2+(k_1\theta_1+k_2\theta_2)\frac{4}{3}\mu_2},\label{M2}
\end{eqnarray}
\begin{eqnarray}
D&=&\Delta T(\frac{3k_1k_2(h_2-h_1)}{k_2-k_1}),\label{D}
\end{eqnarray}
and 
\begin{eqnarray}
F&=&D(1-\frac{1}{\frac{L_1+M_2}{2}}).
\label{F}
\end{eqnarray}

For elastically well--ordered materials, $k_1>k_2$
one has $L_1>1>L_2$, and $M_1>1>M_2$; for the non well--ordered case $k_2>k_1$, one has $L_2>1>L_1$ and $M_2>1>M_1$.

The lower bounds are shown to be obtained by the local fields inside the coated sphere assemblages introduced in \cite{hscoated}. 
The lower bounds  presented here include the effects of thermal stresses due to thermal loads and reduce to the optimal bounds reported in  \cite{Lipstress} when $\Delta T=0$.

\section{Lower bounds on local stress for elastically well-ordered thermoelastic composite media}
\label{sec-optstresshydrowell}
In this section  it is
assumed that the materials inside the heterogeneous medium are elastically well-ordered, i.e., $\mu_1>\mu_2$ and $\kappa_1>\kappa_2$.
We present lower bounds that are optimal for the full range of imposed hydrostatic stresses, i.e., $-\infty<\sigma_0<\infty$ as well as for unrestricted choices of $\Delta T$. The configurations that attain the bounds are given by the coated sphere assemblages \cite{hscoated}. To fix ideas we describe the coated sphere assemblage made from a core of material one with a coating of material two.
We first fill the cube $Q$ with an assemblage of spheres with sizes ranging down to the infinitesimal. Inside each sphere
one places a smaller concentric  sphere  filled with ``core'' material
one and the surrounding coating is filled with material two. The volume fractions of material one and two are taken to be the same for all of the coated spheres. 

In what follows we list the lower bounds for the well ordered case. These bounds are derived in  Section 5.  Their optimality follows from explicit formulas for the moments of the local fields inside the coated sphere assemblage, these are discussed and presented  in Section 5.
The first set of bounds apply to all moments $\langle|\bbP^H\sigma|^p\rangle_2^{1/p}$ for $1< p\leq\infty$.
We suppose that $h_2>h_1$ fix $\Delta T$ and list the bounds as a function of the imposed macroscopic stress $\sigma_0$. The bounds are displayed in the following table where the optimal microstructures are given by the coated spheres construction. The coating and core phase of the optimal configuration is listed in the table below.

  $$
  \begin{array}{|l |l|p{1.35in}|}
   \hline
   Range & Lower Bound & Optimal microstructure\\ \hline
    -\infty<\sigma_0\leq D & \langle|\bbP^H\sigma|^p\rangle_2^{1/p}\geq 
    \sqrt{3}[(D-\sigma_0)L_2-D] &  Core material 2 and  coating  material 1 \\ \hline
   D\leq \sigma_0\leq D(1-\frac{1}{M_2}) & \langle|\bbP^H\sigma|^p\rangle_2^{1/p}\geq \sqrt{3}[(D-\sigma_0)M_2-D]& Core material 1 and coating  material 2\\ \hline
   D(1-\frac{1}{M_2})< \sigma_0< D(1-\frac{1}{L_2}) & \langle|\bbP^H\sigma|^p\rangle_2^{1/p}\geq 0& Optimality undetermined \\ \hline
   D(1-\frac{1}{L_2})\leq \sigma_0<\infty & \langle|\bbP^H\sigma|^p\rangle_2^{1/p}\geq \sqrt{3}[(\sigma_0-D)L_2+D]& Core material 2 and coating material 1\\ \hline
   \end{array}
   $$

Next we suppose that $h_1>h_2$ and  present optimal lower bounds on $\langle|\bbP^H\sigma|^p\rangle_2^{1/p}$, for $1< p\leq\infty$. The bounds and associated optimal microstructures are given in the following table.
$$
\begin{array}{|l |l|p{1.35in}|}
   \hline
   Range & Lower Bound & Optimal microstructure\\ \hline
   -\infty<\sigma_0\leq D(1-\frac{1}{L_2})& \langle|\bbP^H\sigma|^p\rangle_2^{1/p}\geq \sqrt{3}[(D-\sigma_0)L_2-D] & Core material 2 and coating material 1 \\ \hline
   D(1-\frac{1}{L_2})< \sigma_0< D(1-\frac{1}{M_2}) & \langle|\bbP^H\sigma|^p\rangle_2^{1/p}\geq 0 & Optimality undetermined\\ \hline
   D(1-\frac{1}{M_2})\leq\sigma_0\leq D &\langle|\bbP^H\sigma|^p\rangle_2^{1/p}\geq \sqrt{3}[(\sigma_0-D)M_2+D]& Core material 1 and coating material 2\\ \hline
   D\leq \sigma_0 <\infty& \langle|\bbP^H\sigma|^p\rangle_2^{1/p}\geq \sqrt{3}[(\sigma_0-D)L_2+D]& Core material 2 and coating material 1\\ \hline
\end{array}
$$

Lower bounds and the associated optimal microstructures for all moments $\langle|\bbP^H\sigma|^p\rangle_1^{1/p}$ for $1< p\leq\infty$ for the case $h_2>h_1$ is given in the following table.
$$
\begin{array}{|l |l|p{1.35in}|}
   \hline
   Range & Lower Bound & Optimal microstructure\\ \hline
   -\infty<\sigma_0\leq D & \langle|\bbP^H\sigma|^p\rangle_1^{1/p}\geq\sqrt{3}[(D-\sigma_0)L_1-D] & Core material 1 and coating material 2 \\ \hline
   D\leq \sigma_0\leq D(1-\frac{1}{M_1}) & \langle|\bbP^H\sigma|^p\rangle_1^{1/p}\geq\sqrt{3}[(D-\sigma_0)M_1-D]& Core material 2 and coating material 1\\ \hline
   D(1-\frac{1}{M_1})< \sigma_0< D(1-\frac{1}{L_1}) & \langle|\bbP^H\sigma|^p\rangle_1^{1/p}\geq 0& Optimality undetermined\\ \hline
   D(1-\frac{1}{L_1})\leq \sigma_0<\infty & \langle|\bbP^H\sigma|^p\rangle_1^{1/p}\geq\sqrt{3}[(\sigma_0-D)L_1+D]& Core material 1 and coating material 2\\ \hline
   
  \end{array}
  $$
Lower bounds and the associated optimal microstructures for all moments $\langle|\bbP^H\sigma|^p\rangle_1^{1/p}$ for $1< p\leq \infty$ for the case $h_1>h_2$ is given in the following table. 

$$
\begin{array}{|l |l|p{1.35in}|}
   \hline
   Range & Lower Bound & Optimal microstructure\\ \hline
   -\infty<\sigma_0\leq D(1-\frac{1}{L_1}) & \langle|\bbP^H\sigma|^p\rangle_1^{1/p}\geq\sqrt{3}[(D-\sigma_0)L_1-D] & Core material 1 and coating material 2 \\ \hline
   D(1-\frac{1}{L_1})< \sigma_0< D(1-\frac{1}{M_1}) & \langle|\bbP^H\sigma|^p\rangle_1^{1/p}\geq 0 & Optimality undetermined\\ \hline
   D(1-\frac{1}{M_1})\leq\sigma_0\leq D & \langle|\bbP^H\sigma|^p\rangle_1^{1/p}\geq\sqrt{3}[(\sigma_0-D)M_1+D]& Core material 2 and coating material 1\\ \hline
   D\leq \sigma_0<\infty & \langle|\bbP^H\sigma|^p\rangle_1^{1/p}\geq\sqrt{3}[(\sigma_0-D)L_1+D]& Core  material 1 and coating  material 2\\ \hline
\end{array}
$$

Next we display lower bounds on $\max_{\x \hbox{\tiny in Q}}\left\{|\bbP^H\sigma(\x)|\right\}$.
We start with the case $h_2>h_1$ and the lower bounds and optimal geometries are given in the following table.
The phase in which the maximum is attained is denoted with an asterisk.
$$
\begin{array}{|l |l|p{1.35in}|}
   \hline
   Range & Lower Bound & Optimal microstructure\\ \hline
   -\infty<\sigma_0\leq D & \max_{\x \hbox{\tiny in Q}}\left\{|\bbP^H\sigma(\x)|\right\}\geq\sqrt{3}[(D-\sigma_0)L_1-D] & Core material $1^*$ and coating  material 2 \\ \hline
   D\leq \sigma_0\leq F & \max_{\x \hbox{\tiny in Q}}\left\{|\bbP^H\sigma(\x)|\right\}\geq\sqrt{3}[(D-\sigma_0)M_2-D]& Core material 1 and coating material $2^*$\\ \hline
   F\leq \sigma_0<\infty & \max_{\x \hbox{\tiny in Q}}\left\{|\bbP^H\sigma(\x)|\right\}\geq\sqrt{3}[(\sigma_0-D)L_1+D]& Core material $1^*$ and coating material 2\\ \hline
  \end{array}
$$
Lower bounds for and optimal microgepmetries for the case $h_1>h_2$ are given in the following table.
The phase in which the maximum is attained is denoted with an asterisk.
$$
\begin{array}{|l|l|p{1.35in}|}
   \hline
   Range & Lower Bound & Optimal microstructure\\ \hline
   -\infty<\sigma_0\leq F & \max_{\x \hbox{\tiny in Q}}\left\{|\bbP^H\sigma(\x)|\right\}\geq\sqrt{3}[(D-\sigma_0)L_1-D] & Core material $1^*$ and coating material 2 \\ \hline
   F\leq \sigma_0\leq D & \max_{\x \hbox{\tiny in Q}}\left\{|\bbP^H\sigma(\x)|\right\}\geq\sqrt{3}[(\sigma_0-D)M_2+D]& Core  material 1 and coating material $2^*$\\ \hline
   D\leq \sigma_0 <\infty & \max_{\x \hbox{\tiny in Q}}\left\{|\bbP^H\sigma(\x)|\right\}\geq\sqrt{3}[(\sigma_0-D)L_1+D]& Core material $1^*$ and coating material 2\\ \hline
\end{array}
$$

\section{Lower bounds on local stress for non well-ordered thermoelastic composite media}
\label{sec-optstresshydrononwell}

In this section  it is
assumed that the materials inside the heterogeneous medium are elastically non well-ordered, i.e., $\mu_1>\mu_2$ and $\kappa_2>\kappa_1$.
We fix $\Delta T$ and present lower bounds that are optimal for the full range of imposed hydrostatic stresses, i.e., $-\infty<\sigma_0<\infty$. The configurations that attain the bounds for the non well-ordered case are also given by the coated sphere assemblages \cite{hscoated}. 

In what follows we list the lower bounds for the non well-ordered case. These bounds are derived in  Section 5.  Their optimality follows from explicit formulas for the moments of the local fields inside the coated sphere assemblage, these are discussed and presented  in Section 5.
The first set of bounds apply to all moments $\langle|\bbP^H\sigma|^p\rangle_2^{1/p}$ for $1<p\leq\infty$.
We suppose that $h_2>h_1$ and list the bounds as a function of the imposed macroscopic stress $\sigma_0$. The bounds are displayed in the following table where the optimal microstructures are given by the coated spheres construction. The coating and core phase of the optimal coated sphere configuration is listed in the table below.

  $$
  \begin{array}{|l |l|p{1.35in}|}
   \hline
   Range & Lower Bound & Optimal microstructure\\ \hline
    -\infty<\sigma_0\leq D(1-\frac{1}{M_2}) & \langle|\bbP^H\sigma|^p\rangle_2^{1/p}\geq 
    \sqrt{3}[(D-\sigma_0)M_2-D] &  Core material 1 and  coating  material 2 \\ \hline
   D(1-\frac{1}{M_2})\leq \sigma_0\leq D(1-\frac{1}{L_2}) & \langle|\bbP^H\sigma|^p\rangle_2^{1/p}\geq 0 & Optimality undetermined\\ \hline
   D(1-\frac{1}{L_2})< \sigma_0< D & \langle|\bbP^H\sigma|^p\rangle_2^{1/p}\geq \sqrt{3}[(\sigma_0-D)L_2+D]& Core material 2 and coating  material 1 \\ \hline
   D\leq \sigma_0<\infty & \langle|\bbP^H\sigma|^p\rangle_2^{1/p}\geq \sqrt{3}[(\sigma_0-D)M_2+D]& Core material 1 and coating material 2\\ \hline
   \end{array}
   $$

Next we suppose that $h_1>h_2$ and  present optimal lower bounds on $\langle|\bbP^H\sigma|^p\rangle_2^{1/p}$, for $1< p\leq\infty$. The bounds and associated optimal microstructures are given in the following table.
$$
\begin{array}{|l |l|p{1.35in}|}
   \hline
   Range & Lower Bound & Optimal microstructure\\ \hline
   -\infty<\sigma_0\leq D& \langle|\bbP^H\sigma|^p\rangle_2^{1/p}\geq \sqrt{3}[(D-\sigma_0)M_2-D] & Core material 1 and coating material 2 \\ \hline
   D< \sigma_0< D(1-\frac{1}{L_2}) & \langle|\bbP^H\sigma|^p\rangle_2^{1/p}\geq \sqrt{3}[(D-\sigma_0)L_2-D] & Core material 2 and coating material 1\\ \hline
   D(1-\frac{1}{L_2})\leq\sigma_0\leq D(1-\frac{1}{M_2}) &\langle|\bbP^H\sigma|^p\rangle_2^{1/p}\geq 0& Optimality undetermined\\ \hline
   D(1-\frac{1}{M_2})\leq \sigma_0 <\infty& \langle|\bbP^H\sigma|^p\rangle_2^{1/p}\geq \sqrt{3}[(\sigma_0-D)M_2+D]& Core material 1 and coating material 2\\ \hline
\end{array}
$$

Lower bounds and the associated optimal microstructures for all moments $\langle|\bbP^H\sigma|^p\rangle_1^{1/p}$ for $1< p\leq\infty$ for the case $h_2>h_1$ is given in the following table.
$$
\begin{array}{|l |l|p{1.35in}|}
   \hline
   Range & Lower Bound & Optimal microstructure\\ \hline
   -\infty<\sigma_0\leq D(1-\frac{1}{M_1}) & \langle|\bbP^H\sigma|^p\rangle_1^{1/p}\geq\sqrt{3}[(D-\sigma_0)M_1-D] & Core material 2 and coating material 1 \\ \hline
   D(1-\frac{1}{M_1})\leq \sigma_0\leq D(1-\frac{1}{L_1}) & \langle|\bbP^H\sigma|^p\rangle_1^{1/p}\geq 0& Optimality undetermined\\ \hline
   D(1-\frac{1}{L_1})< \sigma_0< D & \langle|\bbP^H\sigma|^p\rangle_1^{1/p}\geq \sqrt{3}[(\sigma_0-D)L_1+D]& Core material 1 and coating material 2\\ \hline
   D\leq \sigma_0<\infty & \langle|\bbP^H\sigma|^p\rangle_1^{1/p}\geq\sqrt{3}[(\sigma_0-D)M_1+D]& Core material 2 and coating material 1\\ \hline
   
  \end{array}
  $$
Lower bounds and the associated optimal microstructures for all moments $\langle|\bbP^H\sigma|^p\rangle_1^{1/p}$ for $1< p\leq \infty$ for the case $h_1>h_2$ is given in the following table.

$$
\begin{array}{|l |l|p{1.35in}|}
   \hline
   Range & Lower Bound & Optimal microstructure\\ \hline
   -\infty<\sigma_0\leq D & \langle|\bbP^H\sigma|^p\rangle_1^{1/p}\geq\sqrt{3}[(D-\sigma_0)M_1-D] & Core material 2 and coating material 1 \\ \hline
   D< \sigma_0< D(1-\frac{1}{L_1}) & \langle|\bbP^H\sigma|^p\rangle_1^{1/p}\geq \sqrt{3}[(D-\sigma_0)L_1-D] & Core material 1 and coating material 2\\ \hline
   D(1-\frac{1}{L_1})\leq\sigma_0\leq D(1-\frac{1}{M_1}) & \langle|\bbP^H\sigma|^p\rangle_1^{1/p}\geq 0& Optimality undetermined\\ \hline
   D(1-\frac{1}{M_1})\leq \sigma_0<\infty & \langle|\bbP^H\sigma|^p\rangle_1^{1/p}\geq\sqrt{3}[(\sigma_0-D)M_1+D]& Core  material 2 and coating  material 1\\ \hline
\end{array}
$$

Next we display lower bounds on $\max_{\x \hbox{\tiny in Q}}\left\{|\bbP^H\sigma(\x)|\right\}$.
We start with the case $h_2>h_1$ and the lower bounds and optimal geometries are given in the following table.
The phase in which the maximum is attained is denoted by an asterisk.
$$ 
\begin{array}{|l |l|p{1.35in}|}
   \hline
   Range & Lower Bound & Optimal microstructure\\ \hline
   -\infty<\sigma_0\leq F & \max_{\x \hbox{\tiny in Q}}\left\{|\bbP^H\sigma(\x)|\right\}\geq\sqrt{3}[(D-\sigma_0)M_2-D] & Core material 1 and coating  material $2^\ast$ \\ \hline
   F\leq \sigma_0\leq D & \max_{\x \hbox{\tiny in Q}}\left\{|\bbP^H\sigma(\x)|\right\}\geq\sqrt{3}[(\sigma_0-D)L_1+D]& Core material $1^\ast$ and coating material 2\\ \hline
   D\leq \sigma_0<\infty & \max_{\x \hbox{\tiny in Q}}\left\{|\bbP^H\sigma(\x)|\right\}\geq\sqrt{3}[(\sigma_0-D)M_2+D]& Core material 1 and coating material $2^\ast$\\ \hline
  \end{array}
$$
Lower bounds for and optimal microgepmetries for the case $h_1>h_2$ are given in the following table.
The phase in which the maximum is attained is denoted by an asterisk.
$$
\begin{array}{|l |l|p{1.35in}|}
   \hline
   Range & Lower Bound & Optimal microstructure\\ \hline
   -\infty<\sigma_0\leq D & \max_{\x \hbox{\tiny in Q}}\left\{|\bbP^H\sigma(\x)|\right\}\geq\sqrt{3}[(D-\sigma_0)M_2-D] & Core material 1 and coating material $2^\ast$ \\ \hline
   D\leq \sigma_0\leq F & \max_{\x \hbox{\tiny in Q}}\left\{|\bbP^H\sigma(\x)|\right\}\geq\sqrt{3}[(D-\sigma_0)L_1-D]& Core  material $1^\ast$ and coating material 2\\ \hline
   F\leq \sigma_0 <\infty & \max_{\x \hbox{\tiny in Q}}\left\{|\bbP^H\sigma(\x)|\right\}\geq\sqrt{3}[(\sigma_0-D)M_2+D]& Core material 1 and coating material $2^\ast$\\ \hline
\end{array}
$$

\bigskip
\section{Derivation of the lower bounds on $<|\mathbb{P}_H\sigma|^p>_i^\frac{1}{p}$}
\label{lower bound 2}

In this section we outline the methodology for proving optimal lower bounds.
The bounds are derived using duality relations.  We use the following duality relation posed over the space of square integrable symmetric matrix fields $\eta$ that holds for $p>1$ given by
\begin{eqnarray}
\frac{1}{p}\langle|\bbP^H\sigma|^p\rangle_i=\sup_{\eta}\left\{\langle\bbP^H\sigma: \eta\rangle_i-\frac{1}{p'}\langle|\eta|^{p'}\rangle_i \right\},\hbox{  for $i=1,2$,
\label{vec1}}
\end{eqnarray}
where $p'$ is the conjugate exponent to $p$ given by $p'=\frac{p}{p-1}$.
This relation follows immediately from standard duality relations see, \cite{Dacorogna}.
Restricting $\eta$ to the set of all constant matrices and taking the supremum delivers the basic bounds:
\begin{eqnarray}
\langle|\bbP^H\sigma|^p\rangle_i&\geq &|\langle\bbP^H\sigma\rangle_i|^p
,\hbox{ $p>1$ and for $i=1,2$}.
\label{vecbound2}
\end{eqnarray}
We point out that equality holds in \eqref{vecbound2} if and only if $\bbP^H\sigma$ is identically constant inside the  $i^{th}$ material.
In what follows we outline the method for obtaining bounds on the moments in material two noting that the identical procedure delivers bounds on the moments in  material one. We introduce the indicator function of material one $\chi_1$ taking the value 1 in material one and zero outside. The indicator function corresponding to material two is denoted by $\chi_2$ and $\chi_2=1-\chi_1$.
To proceed we rewrite the right hand side of \eqref{vecbound2} in terms of the effective elastic properties and thermal expansion coefficient.
To do this we use the following identity given by
\begin{align}\label{5.1}
 tr<\chi_2\sigma>=\frac{3k_2}{k_2-k_1}(\sigma_0-k_1\sigma_0(C^e)^{-1}I:I+k_1\Delta T (C^e)^{-1}H^e:I+k_1\Delta T<\lambda>:I).
\end{align}
\par
This identity is obtained in the following way. 
Taking averages on both sides of  \eqref{constitutive}  gives $<\sigma>=\langle C(x)(\epsilon (\textbf{u}(x))-\Delta T\lambda(x))\rangle$. On writing $C(\x)=C_1+(C_2-C_1)\chi_2(\x)$ we see that
\begin{align}\label{5.2}
 <\sigma>=C_1(\bar{\epsilon}-<\lambda>)+(C_2-C_1)C_2^{-1}<\chi_2\sigma>.
\end{align}

From \eqref{effelast} we see that $\bar{\epsilon}=(C^e)^{-1}(<\sigma>-\Delta T H^e)$ and  for $\langle\sigma\rangle=\sigma_0 I$ we obtain
\begin{align}\label{5.3}
<\chi_2\sigma>=C_2(C_2-C_1)^{-1}(\sigma_0I-C_1((C^e)^{-1}\sigma_0 I-\Delta T (C^e)^{-1}H^e-\Delta T<\lambda>)).
\end{align}
The identity (\ref{5.1}) now follows by applying the hydrostatic projection $\mathbb{P}_H$ to both sides of (\ref{5.3}).
\par
We now derive the lower bound. Applying the basic bound \eqref{vecbound2} to $
\langle|\mathbb{P}_H\sigma|^p\rangle_2$ and (\ref{5.1}) gives
\begin{eqnarray}
&&\langle|\mathbb{P}_H\sigma|^p\rangle_2\nonumber\\
&&\geq |\langle\mathbb{P}_H\sigma \rangle_2|^p=\left(\frac{tr(\langle\chi_2\sigma\rangle)}{\sqrt{3}\theta_2}\right)^p\nonumber\\
&&=3^{p/2}\theta_2^{-p}|\frac{k_2}{k_2-k_1}|^p|k_1\sigma_0\left(\frac{1}{k_1}-(C^e)^{-1}I:I\right)+k_1\Delta T(C^e)^{-1}H^eI:I+k_1\Delta T \langle\lambda\rangle:I|^p.\label{5.4}
\end{eqnarray} 
%
%
We note that equality holds in (\ref{5.4}) when $\mathbb{P}_H\sigma$ is constant inside material two.

We now employ an exact relation that relates  the contraction $(C^e)^{-1}H^e:I$ involving the effective thermal stress tensor $H^e$ to the quantity $(C^e)^{-1}I:I$. The exact relation used here is given by
\begin{eqnarray}
(C^e)^{-1}H^e:I&=&\frac{3(h_2-h_1)(C^e)^{-1}I:I+3(\frac{h_1}{k_2}-\frac{h_2}{k_1})}{\frac{1}{k_1}-\frac{1}{k_2}}.
\label{exact}
\end{eqnarray}
This exact relation is a direct consequence of the exact relation developed by  \cite{HashinRosen} for the effective thermal expansion tensor $\mathbf{\alpha}^e=-(C^e)^{-1}H^e$.

Substitution of \eqref{exact} into \eqref{5.4} and algebraic manipulation gives
\begin{eqnarray}
\langle|\mathbb{P}_H\sigma|^p\rangle_2^{1/p}\geq \sqrt{3}\left|(\sigma_0-D)X+ D\right|,
\label{finalbd}
\end{eqnarray}
where
\begin{eqnarray}
X=\frac{\theta_2^{-1}}{k_1^{-1}-k_2^{-1}}\left(\frac{1}{k_1}-(C^e)^{-1}I:I\right).
\label{x}
\end{eqnarray}
As before we point out that equality holds in (\ref{finalbd}) when $\mathbb{P}_H\sigma$ is identically constant inside material two.

Identical arguments deliver the lower bound on the moments over material two given by
\begin{eqnarray}
\langle|\mathbb{P}_H\sigma|^p\rangle_1^{1/p}\geq \sqrt{3}\left|(\sigma_0-D)Y+ D\right|,
\label{finalbd2}
\end{eqnarray}
where
\begin{eqnarray}
Y=\frac{\theta_1^{-1}}{k_2^{-1}-k_1^{-1}}\left(\frac{1}{k_2}-(C^e)^{-1}I:I\right),
\label{y}
\end{eqnarray}
where \eqref{finalbd2} holds with equality when $\bbP^H\sigma$ is a constant in material one.

The variables $X$ and $Y$ are constrained to lie within intervals set by bounds on the contraction
$(C^e)^{-1}I:I$. These bounds follow immediately from the work of Kantor and Bergman \cite{KB} and are given by
\begin{eqnarray}
(K_{HS}^+)^{-1}\leq (C^e)^{-1}I:I \leq (K_{HS}^-)^{-1},
\label{HS}
\end{eqnarray}
where $K_{HS}^-$ and $K_{HS}^+$ are the Hashin and Shtrikman  bulk modulus bounds \cite{hashin} given by 
\begin{eqnarray}
K_{HS}^+=k_1\theta_1+k_2\theta_2-(\frac{\theta_1\theta_2(k_2-k_1)^2}{k_1\theta_2+k_2\theta_1+\frac{4}{3}\mu_1})
\label{hs+}
\end{eqnarray} 
and
\begin{eqnarray} K_{HS}^-=k_1\theta_1+k_2\theta_2-(\frac{\theta_1\theta_2(k_2-k_1)^2}{k_1\theta_2+k_2\theta_1+\frac{4}{3}\mu_2}).
\label{hs-}
\end{eqnarray}
These bounds hold both for elastically well-ordered materials and elastically non well-ordered materials.
When the materials are well--ordered \eqref{HS} implies that $X$ and $Y$ lie in the intervals
\begin{eqnarray}
L_2\leq & X&\leq  M_2,\label{intervalwell1}\\
L_1\leq & Y& \leq  M_1.
\label{intervalwell2}
\end{eqnarray}
while for non well--ordered materials 
\begin{eqnarray}
M_2\leq & X&\leq L_2,\label{intervalnonwell1}\\
M_1\leq & Y& \leq  L_1.
\label{intervalnonwell2}
\end{eqnarray}

A straightforward calculation in Section \ref{Newsolutions} shows that the hydrostatic component of the local stress is constant inside each phase of the coated sphere construction.  Hence  \eqref{finalbd} and \eqref{finalbd2} hold with equality for the coated spheres construction and we obtain explicit formulas for the moments of the hydrostatic component of the local stresses for these composites. For coated spheres with core phase 1 and coating phase 2, $(C^e)^{-1}I:I=(K_{HS}^-)^{-1}$ and
substitution of  \eqref{hs-} into \eqref{x} and \eqref{y} together with  \eqref{finalbd} and  \eqref{finalbd2}  shows that the moments are given by
\begin{eqnarray}
\langle|\mathbb{P}_H\sigma|^p\rangle_2^{1/p}=\sqrt{3}\left|(\sigma_0-D)M_2+ D\right|,\hbox{ and}
\label{finalbdhs+}\\
\langle|\mathbb{P}_H\sigma|^p\rangle_1^{1/p}= \sqrt{3}\left|(\sigma_0-D)L_1+ D\right|.
\label{finalbd2hs+}
\end{eqnarray}
For coated spheres with core phase 2 and coating phase 1, $(C^e)^{-1}I:I=(K_{HS}^+)^{-1}$ and
substitution of  \eqref{hs+} into \eqref{x} and \eqref{y} together with \eqref{finalbd} and  \eqref{finalbd2} shows that the moments are given by
\begin{eqnarray}
\langle|\mathbb{P}_H\sigma|^p\rangle_2^{1/p}=\sqrt{3}\left|(\sigma_0-D)L_2+ D\right|,\hbox{ and}
\label{finalbdhs-}\\
\langle|\mathbb{P}_H\sigma|^p\rangle_1^{1/p}= \sqrt{3}\left|(\sigma_0-D)M_1+ D\right|.
\label{finalbd2hs-}
\end{eqnarray}

We collect results and state the lower bounds and indicate when they are optimal.

\begin{theorem}{\rm Bounds for well-ordered composites, $k_1>k_2$.}
\label{optlowerboundswell}
For $1< p \leq\infty$, any choice of $\Delta T$,  and $-\infty<\sigma_0<\infty$ the lower bounds are given by the following formulas.
\begin{eqnarray}
\langle|\mathbb{P}_H\sigma|^p\rangle_2^{1/p}\geq \min_{L_2\leq X\leq  M_2}\left\{\sqrt{3}\left|(\sigma_0-D)X+ D\right|\right\},
\label{finalbdmin2well}
\end{eqnarray}
and when the minimum is realized for $X=L_2$ the bound is attained by the fields inside the core phase of a coated sphere construction with core material 2 and coating 1;  when the minimum is realized for $X=M_2$ the bound is attained by the fields inside the coating phase of a coated sphere construction with core material 1 and coating 2. 

\begin{eqnarray}
\langle|\mathbb{P}_H\sigma|^p\rangle_1^{1/p}\geq \min_{L_1\leq Y\leq  M_1}\left\{\sqrt{3}\left|(\sigma_0-D)Y+ D\right|\right\},
\label{finalbdmin1well}
\end{eqnarray}
and when the minimum is realized for $Y=L_1$ the bound is attained by the fields inside the core phase of a coated sphere construction with core material 1 and coating 2 and  when the minimum is realized for $Y=M_1$ the bound is attained by the fields inside the coating phase of a coated sphere construction with core material 2 and coating 1.
\end{theorem}

\begin{theorem}{\rm Bounds for non well-ordered composites, $k_2>k_1$.}
\label{optlowerboundsnonwell}
For $1< p \leq\infty$, any choice of $\Delta T$, and $-\infty<\sigma_0<\infty$ the lower bounds are given by the following formulas.
\begin{eqnarray}
\langle|\mathbb{P}_H\sigma|^p\rangle_2^{1/p}\geq \min_{M_2\leq X\leq  L_2}\left\{\sqrt{3}\left|(\sigma_0-D)X+ D\right|\right\},
\label{finalbdmin2nonwell}
\end{eqnarray}
and when the minimum is realized for $X=L_2$ the bound is attained by the fields inside the core phase of a coated sphere construction with core material 2 and coating 1;  when the minimum is realized for $X=M_2$ the bound is attained by the fields inside the coating phase of a coated sphere construction with core material 1 and coating 2. 

\begin{eqnarray}
\langle|\mathbb{P}_H\sigma|^p\rangle_1^{1/p}\geq \min_{M_1\leq Y\leq  L_1}\left\{\sqrt{3}\left|(\sigma_0-D)Y+ D\right|\right\},
\label{finalbdmin1nonwell}
\end{eqnarray}
and when the minimum is realized for $Y=L_1$ the bound is attained by the fields inside the core phase of a coated sphere construction with core material 1 and coating 2 and  when the minimum is realized for $Y=M_1$ the bound is attained by the fields inside the coating phase of a coated sphere construction with core material 2 and coating 1.
\end{theorem}

These bounds are stated explicitly in the first four tables of Sections 3 and 4.

We conclude by outlining the steps behind the derivation of the lower bounds on the maximum values of the local fields inside thermally stressed composites. For the well--ordered case we use the simple lower bound given by
\begin{eqnarray}
\max_{ \x\hbox{ \tiny in $Q$}}\left\{\right|\bbP^H\sigma(\x)\}\geq\max\left\{A,B\right\}.
\label{linftyab}
\end{eqnarray}
where
\begin{eqnarray}
A&=&\min_{L_2\leq X\leq  M_2}\left\{\sqrt{3}\left|(\sigma_0-D)X+ D\right|\right\},\nonumber\\
B&=&\min_{L_1\leq Y\leq  M_1}\left\{\sqrt{3}\left|(\sigma_0-D)Y+ D\right|\right\}.
\label{ab}
\end{eqnarray}For the non well--ordered case we use
\begin{eqnarray}
\max_{\x\hbox{ \tiny  in $Q$}}\left\{\right|\bbP^H\sigma(\x)\}\geq\max\left\{C,D\right\}.
\label{linftycd}
\end{eqnarray}
where
\begin{eqnarray}
C&=&\min_{M_2\leq X\leq  L_2}\left\{\sqrt{3}\left|(\sigma_0-D)X+ D\right|\right\},\nonumber\\
D&=&\min_{M_1\leq Y\leq  L_1}\left\{\sqrt{3}\left|(\sigma_0-D)Y+ D\right|\right\}.
\label{cd}
\end{eqnarray}

The bounds given in the last two tables presented in Sections 3 and 4 follow from  straight forward but tedious calculation of the explicit formulas corresponding to \eqref{linftyab} and \eqref{linftycd}.
A delicate but straight forward computation shows that these lower bounds are attained by the fields inside the coated sphere assemblage.

\bigskip
\bigskip
\section{Local stress and strain fields inside  thermally stressed  coated sphere geometries}
\label{Newsolutions}


In this section we summarize the properties of local fields inside the coated sphere assemblage in the presence of thermal stress due to a mismatch in the coefficients of thermal expansion. From linearity the local stress can be split into the sum of two components; one component arising from imposed mechanical stress and a second component associated with thermal stress. It is known that the local stress due to an imposed hydrostatic stress has constant hydrostatic part inside each phase, this follows from explicit solution see for example \cite{Milton}.
Here we display the  explicit solution for the local stress due to mismatch in the coefficients of thermal expansion and show that it has a constant hydrostatic component inside each phase. From this we conclude that the total local stress inside the coated sphere assemblage has a constant hydrostatic component inside each phase. 

We solve for the stress inside a prototypical coated sphere composed of a spherical core of material two with  radius $a$, surrounded by a concentric shell of material one with an outer radius $b$. The ratio $(a/b)^3$ is fixed and equal to the inclusion volume fraction $\theta_2$. Here the the coefficients of thermal expansion for the core and coating are given by $h_1$ and $h_2$ respectively. The local elastic displacement $\tilde{\varphi}$ satisfies the equations of elastic equilibrium are given by:

\begin{displaymath}
\left\{\begin{array}{ll} div(C_2(\epsilon (\tilde{\varphi})-h_2I))=0 & 0<r<a ,\\
div(C_1(\epsilon (\tilde{\varphi})-h_1I))=0 & a<r<b,\\
C_1(\epsilon (\tilde{\varphi})-h_1I) \textbf{n}|_1=C_2(\epsilon (\tilde{\varphi})-h_2I)\textbf{n}|_2 & \text{continuity of traction at } r=a,\\
\tilde{\varphi} \quad\text{is continuous}& \text{on $0<r<b$},\\
\tilde{\varphi}=0  &\text{on the boundary $r=b$}.\end{array}\right.
\end{displaymath} 
\par
We assume a general form of the solution given by
\begin{equation*}
\tilde{\varphi} =
\begin{cases}
C\textbf{r} & 0<r<a,\\
A\textbf{r}+B\frac{\textbf{n}}{r^2} & a<r<b,\\
0 & r\geq b,
\end{cases}
\end{equation*}
where $r=|\textbf{r}|, \textbf{n}=\textbf{r}/r$ and $A, B, C$ are unknowns. The corresponding strain field $\epsilon(\tilde{\varphi})$ is given by
\begin{equation}
(\epsilon (\tilde{\varphi}))_{ij} =\frac{1}{2}(\tilde{\varphi}_{i,j}+\tilde{\varphi}_{j,i})=
\begin{cases}
C\delta_{ij} & 0<r<a,\\
A\delta_{ij}+\frac{B}{r^3}(\delta_{ij}-3n_in_j) & a<r<b,\\
0 & r\geq b.
\end{cases}
\label{hydroconst}
\end{equation}

On applying the continuity of displacement and the traction at the interface we find that
$$A=\frac{3\theta_2(k_1h_1-k_2h_2)}{3k_1\theta_2+4\mu_1+3k_2(1-\theta_2)} , $$
$$B=\frac{-3a^3(k_1h_1-k_2h_2)}{3k_1\theta_2+4\mu_1+3k_2(1-\theta_2)} , $$
$$C=\frac{-3(1-\theta_2)(k_1h_1-k_2h_2)}{3k_1\theta_2+4\mu_1+3k_2(1-\theta_2)}.$$
\\
\par
Computation of the radical component of the stress at $r=b$ gives
\begin{eqnarray}
C_1(\epsilon(\tilde{\varphi})-h_1I) \textbf{n}=H^\ast\textbf{n}
\label{normal}
\end{eqnarray}
where $H^\ast$
\begin{align}\label{2.3}
H^\ast=\frac{3\theta_2(3k_1+4\mu_1)(k_1h_1-k_2h_2)}{3k_1\theta_2+4\mu_1+3k_2(1-\theta_2)}I-3k_1h_1I.
\end{align}
\\
\par
In this way we have constructed a solution $\tilde{\varphi}$ for the elastic field inside  every coated sphere in the assemblage.
We now define $\varphi^p$ on the whole domain $Q$ to be given by $\tilde{\varphi}$ inside each coated sphere and zero outside. It easily follows on integrating by parts using
\eqref{normal} together with the fact that $\varphi^p$ vanishes on the boundary of each coated sphere that $\varphi^p$ is the weak  solution \cite{Gilbarg} of $div(C(\epsilon (\varphi^p-\lambda))=0 $ over the full domain $Q$, i.e.,$$\langle C(\epsilon (\varphi^p)-\lambda):\epsilon(\phi)\rangle=0$$ for every periodic test function $\phi$.
Equation \eqref{hydroconst} implies that hydrostatic component of stress is constant inside each phase. 

Last we show that the effective thermal stress $H^e$ for the coated sphere assemblage is given by $H^\ast$. Inside each coated sphere $S_i$, $i=1,2,\ldots$ we consider $\sigma=C(\epsilon (\varphi^p)-\lambda)$ and integrate by parts and apply \eqref{normal} to find that
\begin{eqnarray}
\int_{S_i}\sigma d\x&=&\int_{\partial S_i}(\sigma\n)\otimes\x ds\\
&=&\int_{\partial S_i}(H^\ast\n)\otimes\x ds=|S_i|H^\ast
\label{intbyparts}
\end{eqnarray}
where $ds$ is an element of surface area on the outer surface of the coated sphere $\partial S_i$
and $|S_i|$ is the volume of $S_i$.
Substitution of \eqref{intbyparts} into \eqref{effectethstress} 
gives the required identity
\begin{eqnarray}
H^e=\langle\sigma\rangle=\sum_{i=1}^\infty \int_{S_i}\sigma d\x=H^\ast.
\label{intbypartsandidentity}
\end{eqnarray}


\begin{thebibliography}{99}

\bibitem{AlaliLipton}
Alali, B. and Lipton, R. ``Optimal lower bounds on local stress inside random media.'' SIAM J. On Applied Math. To Appear electronically, November, 2009.


\bibitem{Dacorogna}
Dacorogna, B., 1989. {\em Direct Methods in the Calculus of Variations}, Springer--Verlag, Berlin.


\bibitem{faraco}
Faraco, D., 2003. ``Milton's conjecture on the regularity of solutions to isotropic equations.''
Annales de L'Institute Henri Poincare (c) Nonlinear Analysis {\bf 20}, pp. 889--909.


\bibitem{Gilbarg}
Gilbarg, D. and Trudinger N., 2001. {\em Elliptic Partial Differential Equations of Second Order,} Springer, Berlin.

\bibitem{GrabovskyandKohn}
Grabovsky, Y. and Kohn, R. V., 1995. ``{Microstructures minimizing the energy of a two phase elastic composite in two space
dimensions. II: The Vigdergauz microstructure.}'' J. Mech. Phys. Solids. {\bf 43}, pp. 949--972.
%



\bibitem{Hash}
Hashin, Z., 1983. ``Analysis of composite materials -- a survey.'' Journal of Applied Mechanics
{\bf 50}, pp. 481--505.



\bibitem{hscoated}
Hashin, Z.,  1962. ``The elastic moduli of heterogeneous materials.'' Journal of Applied Mechanics {\bf 29}, pp. 143--150.

\bibitem{hashin}
Hashin, Z. and Shtrikman, S., 1963. ``A variational approach to the theory of the elastic behaviour of multiphase materials.'' J. Mech. Phys. Solids, {\bf 11}, pp. 127--140.

\bibitem{He}
He, Q.C., 2007. ``Lower bounds in the stress and strain fields inside random two-phase elastic media.''
Acta Mechanica {\bf 188}, pp. 123--137.

\bibitem{Hill}
Hill, R., 1963. ``Elastic properties of reinforced solids: Some theoretical principles.'' J. Mech. Phys. Solids {\bf 11}, pp. 357--372.


\bibitem{KB}
Kantor, Y. and Bergman, D.J. 1984. ``Improved rigorous bounds on the effective elastic moduli of a
composite material.'' J. Mech. Phys. Solids. {\bf 32}, pp. 41--62.



\bibitem{leonetinesi}
Leonetti, F. and Nesi, V., 1997. ``Quasiconformal solutions to certain first order systems and the proof of a conjecture of G.W. Milton.'' J. Math. Pures. Appl. {\bf 76}, pp. 109--124.





\bibitem{Lipelect}
Lipton, R., 2004. ``Optimal lower bounds
on the electric-field concentration in composite media.'' { Journal of
Applied Physics}, {\bf 96}, pp. 2821--2827.

\bibitem{Lipstrain}
Lipton, R., 2006.``Optimal lower bounds on the dilatational strain inside random
two-phase elastic composites subjected to hydrostatic loading.'' {Mechanics of Materials}, {\bf 38}, pp. 833--839.

\bibitem{Lipstress}
Lipton, R., 2005. ``Optimal lower bounds on the hydrostatic stress amplification inside random
two-phase elastic composites.'' {Journal of the Mechanics and Physics of Solids} {\bf 53},
pp. 2471--2481.

\bibitem{lipsiama}
Lipton, R., 2001. ``Optimal inequalities for gradients of solutions of elliptic equations occurring in two-phase heat conductors. SIAM J. Math. Analysis, {\bf 32}, pp. 1081--1093.

\bibitem{Maxwell04}
Maxwell Garnett, J.C., 1904. ``Colours in metal glasses and in metallic films. Philosophical
Transactions of the Royal Society of London {\bf 203}, pp. 385--420.


\bibitem{Milton}
Milton, G.W., 2002. {\em The Theory of Composites}. Cambridge University Press, Cambridge.

\bibitem{SeminalMilton}
Milton, G.W., 1986. ``Modeling the properties of composites by laminates.'' In: Homogenization
and Effective Moduli of Materials and Media. Edited by J. Erickson, D. Kinderleher, R.V. Kohn, and J.L. Lions.
IMA Volumes in Mathematics and Its Applications 1, pp. 150--174. Springer-Verlag, New York.

\bibitem{NNH}
Nemat--Nasser, S. and Hori, M., 1999. {\em Micromechanics: Overall Properties
of Heterogeneous Materials}. Elsevier, Amsterdam.



\bibitem{Rayleigh}
Rayleigh, L., 1892. On the influence of obstacles arranged in rectangular order upon the
properties of a medium. Philosophical Magazine {\bf{34}}, pp. 481--502.

\bibitem{HashinRosen}
Rosen, B. W., and Hashin, Z., 1970. ``Effective thermal expansion coefficients and specific heats of composite materials.'' International Journal of Engineering Science {\bf{8}}, pp. 157--173. 






\bibitem{Wheeler}
Wheeler, L. T., 1993. `` Inhomogeneities of minimum stress concentration.'' Anisotropy and Inhomogeneity in
Elasticity and Plasticity, Y. C. Angel, ed., AMD-Vol. 158. ASME, pp. 1--6.

\bibitem{Willis}
Willis, J. R., 1983. ``The overall elastic response of composite materials.'' J. App. Mech., {\bf 50}, pp. 1202--1209.


\end{thebibliography}
\end{document}